# A Note on the Possibility of Self-Reference in Mathematics

## Arieh Lev

**Abstract.** In this paper we propose an interpretation for self-referential sentences and propositions in a 'meta- model' $N^*$ of ZF. This meta model $N^*$ is considered as an 'informal' model of arithmetic that mathematicians often use when working with number theory. Specifically, we assume that within this meta-model, the axiom system ZF is applied, interpretations for sentences can be offered, and natural language can be used. We show that under the proposed interpretation, some types of self-referential sentences that are considered legitimate in mathematics turn $N^*$ into an inconsistent model, and examine the connection of this result to a certain interpretation of Gödel's first incompleteness theorem. Some general problems which follow from the above discussion are then addressed.

### 1) Introduction

In this paper we examine the significance of sentences that refer to themselves, in an axiomatized theory $T$ that includes the Zermelo-Fraenkel axiom system for set theory (ZF). We also address propositions that constitute an interpretation of these sentences.

Specifically, we address propositions which refer to their truth value or their provability, like the following two sentences:

$H$ : The proposition $H$ is a true proposition.

$G$ : The proposition $G$ is unprovable.

In this paper, we denote the 'standard' model for arithmetic by $N$. We also denote the '*meta-model*' of arithmetic, that mathematicians generally use when working with number theory, by $N^*$. Specifically, we assume that within this meta-model, the axiom system ZF is considered, interpretations for sentences can be offered, and natural language can be used (for instance, we can use concepts like 'proof,' say that a proposition is true or untrue, provable or unprovable, etc.).

Note that the ZF axiom system meets the requirements of Gödel's incompleteness theorems. We also assume that $T$ is a consistent theory (that it does not imply a proof of both a sentence and its negation).

Our motivation for concerning ourselves with the proposition $G$ above is derived, naturally, from Gödel's first incompleteness theorem. In proving this theorem, Gödel presented a sentence $G_T$, where neither the sentence nor its negation are provable in the theory $T$.

On the interpretation of the sentence $G_T$, Gödel writes: 'We therefore have before us a proposition that says about itself that it is not provable' (see, for instance, [H, p.

598]). In other words, we could say that Gödel showed how we could, within the theory $T$, build sentences that can be interpreted as propositions that speak on their own provability. In light of this, it is common today to say that certain propositions that refer to themselves are possible in mathematics.

It is possible in the meta-model $N^*$ to speak about self-reference. We can therefore speak within it of propositions like the aforementioned $G$, and consider implications of the '*meta-interpretation*' of $G_T$ in $N^*$.

We will focus on the significance of the concept of self-reference in the meta-model $N^*$. We will begin (in Sections 2 and 3 below) by clarifying the meaning of the concepts *self-reference* and *circularity* as they are used in this paper. Then we argue (in Sections 4 and 5) that in the meta-model $N^*$, propositions that refer to themselves (in the sense applied to this in our paper, and assuming the existence of 'natural coding' in their semantics, as described in Section 4), leads to a 'conflict' with the ZF's axiom of foundation. This latter point is of great interest to the study of paradoxes and circularity, even without the relation to Gödel's theorem. On the other hand, our discussion can also be related, at least to some extent, to Gödel's first incompleteness theorem. Though it does not detract from the importance of Gödel's first incompleteness theorem, accepting the point of view presented in this paper may impact the conclusions that arise from this theorem on the meta-model $N^*$. This latter issue will be discussed in Sections 6-8.

The paper is designed to be accessible to a wider audience of readers, including those that are not familiar with topics like self-reference, Gödel's theorem, or the *ZF* axiom system. For more detailed information, the reader may consult [S], [F] and [V].

### 2) Strong self-reference

Let us open this section with a clarification of the concept of self-reference, and begin with an example.

Examine the proposition $S$, which says, 'I am a proposition that includes fewer than 15 words.' This is a proposition that, allegedly, refers to itself. Once we count the words in that proposition, we will immediately be convinced of its validity. Though the proposition may be thought to contain self-reference, a closer inspection may convince us that it does not: the proposition does not speak about itself, but about the text that describes it. More specifically, the truth value of this proposition could change based on the manner or the language in which it is phrased.

In contrast, propositions like 'I am a false proposition' or 'It is impossible to prove me' can be said to be self-referential, since they refer to the proposition's semantics. In particular, the meaning of the propositions will not change with changes to the phrasing or language in which they are written.

We must therefore be careful to address propositions that refer to an 'extensional property' (see e.g., [R, p. 186] or [J, p. 79]): a property that does not refer to the structure or manner of the proposition's description, but rather to its semantics. We will not go into an in-depth discussion of the problem of rigorously defining the 'extensional properties' of propositions here. We shall merely determine that the truth value of a proposition, or its provability, are extensional properties, and that propositions that address their own truth value or provability are therefore to be considered self-referential.

In order to simplify matters, we add another condition to the concept of self-reference in this paper: a condition that demands that the **only** claim made by a proposition will be about its own extensional property. Such a proposition will be said to be **strongly self-referential.** Thus, for instance, the sentence, 'I am a false proposition' fulfils that condition, while the sentence, '2+2=4 is true and I am a true proposition' does not.

**Note:** It may be argued that a proposition which refers to its provability, like '$G$ : The proposition $G$ is unprovable' does not refer to itself, since a proof of a proposition refers only to a syntactic object that represents the proposition, and not to the 'content' or 'meaning' of the proposition. However, we shall argue in Chapter 5 that we may consider another interpretation of $G$ in $N^*$, under which $G$ may be viewed as a strongly self-referential proposition. See Chapter 5 for further details.

From this point onward, whenever we argue that a proposition is self-referential, we will be referring to strong self-reference. We shall say that a sentence (in the theory $T$) is self-referential if its meta-interpretation in $N^*$ is a self-referential proposition. We will be applying similar limitations to our discussion of propositions that refer to the extensional properties of *other* propositions: we will assume that the **only** claim made by that proposition is a claim regarding that extensional property.

3) **Vicious circle**

Russel's paradox is based on the possibility that a set is a member of itself. In this case, Russel proves that it is possible to define a set that is both a member and not a member of itself, thereby producing a paradox. For the purposes of discussing such a circular situation, Russel defined the principle of the *vicious circle*. Russel phrased this principle in rather vague terms (see e.g., [Chi, p. 3]), stating that: 'Whatever involves all of a collection must not be one of the collection.'

Mathematicians could not accept a situation in which set theory contains a 'vicious circle' paradox. They therefore looked for a way to address the problem. Russel himself drew up his theory of types (which is not currently popular in mathematics) specifically for this purpose, and many other attempts at a solution have been made. Today, the accepted approach is Zermelo-Fraenkel's axiomatic set theory (ZF, or ZFC if the choice axiom is included). It has no proof, but mathematicians believe that this axiom system prevents the existence of paradoxes in set theory.

Of particular relevance to our topic is the **axiom of foundation** (also known as the axiom of regularity). We will use the following phrasing of the axiom, which is equivalent to the common phrasing for infinite sets (see, e.g., [V, chapter 8, p. 93]):

**Axiom:** There is no ordinary infinite $\in$-decreasing sequence (that is, there is no infinite sequence of the form $\ldots \in A_3 \in A_2 \in A_1$).

Note that the axiom prevents the possibility that a set will be a member of itself. Indeed, if $A \in A$ is true, this would necessarily lead to the infinite sequence $\ldots \in A \in A \in A$. Moreover, it also prevents situations in which a set $B$ is a member of a set $A$ when the set $A$ is a member of $B$, etc.

To avoid vague phrasing (like that employed in Russel's definition of the vicious circle), we clarify that whenever we claim that a set is contaminated by a vicious circle situation, we mean that assuming the existence of this set results in a conclusion that contradicts the axiom of foundation.

We also apply the principle of the vicious circle to propositions. We say that the principle of the vicious circle applies to a proposition if (the meta-interpretation of) $S$ can produce a conclusion that contradicts the axiom of foundation.

### 4) Self-referential propositions

The discussion in this section takes place within the meta-model of number theory, which we have designated above as $N^*$. Recall that in this meta-model we assume the ZF axiomatic set theory. Our goal in this section is to show that, based on certain (natural) assumptions, the existence of self-referential propositions in $N^*$ stands in contradiction to the axiom of foundation.

For the sake of simplicity, we only address propositions the content of which is wholly devoted to claims about the provability (or non-provability) of given propositions, such as the following proposition $H$, which addresses the provability of proposition $K$:

$H$: The proposition $K$ is unprovable.

Let us examine the following definition:

**Definition:** A set $A$ of propositions will be called a *set of type $P$*, if every proposition in $A$ is a proposition wholly devoted to a statement about the provability (or non-provability) of a proposition belonging to $A$.

Even though the definition of the set of type $P$ may seem 'strange' or 'circular,' there do exist sets of type $P$ that are not empty. Thus, for instance, if $A = \{L\}$ when $L$ is the proposition 'the proposition $L$ is not provable,' then $A$ is a non-empty set of type $P$.

We now fix a set $A_p$, which is a non-empty set of type $P$.

Note that we may interpret a proposition of the form '$H$: The proposition $K$ is provable/unprovable' in $A_p$ as a function $f_H$, where $K$ is the only member of its domain, and its value is 1 if $H$ claims that $K$ is provable, and 0 if $H$ claims that $K$ is unprovable. But $K$ itself is a member of $A_p$, and hence may be also interpreted as a function $f_K$. Therefore, we conclude that $H$ may be interpreted as a function (with domain $\{f_K\}$): $f_H : \{f_K\} \to \{0,1\}$.

Now we define the set $F_p$ of functions as follows. For every proposition '$L$: The proposition $M$ is provable' (or, alternatively, '$L$: the proposition $M$ is unprovable') in the set $A_p$, there is a unique function denoted $f_L$, which is defined as follows: $f_L : \{f_M\} \to \{0,1\}$ where $f_L(f_M) = 1$ if $L$ claims that $M$ is provable and $f_L(f_M) = 0$ otherwise. This means that for every function $f_L \in F_p$, there is a single element within its domain $f_M$, which is itself a function belonging to $F_p$, so that the value $f_L(f_M)$ is determined by the proposition $L$: the value is 1 if $L$ claims that $M$ is provable, and 0 if $L$ claims that $M$ is unprovable.

The natural map from $A_P$ to $F_p$, which transfers every proposition $L$ to the function $f_L$, is one to one and onto. Since the set $A_P$ is not empty, the set $F_p$ is not empty either.

Now let us assume the existence of the following self-referential proposition in the set $A_p$:

$H$: $H$ is not provable.

This leads us to conclude that for the function $f_H : \{f_H\} \to \{0,1\}$, we have $f_H(f_H) = 0$, and specifically, we have that $f_H$ is a member of its own domain.

Let us recall that every function $f : A \to \{0,1\}$ can be presented as a set of ordered pairs $(a,b)$ (known as the graph of the function, see [V p. 16]), where $(a,b)$ is a member of the graph of $f$ if and only if $a \in A$ and $b \in \{0,1\}$ is the only member for which $f(a) = b$. We can therefore, as is often done, represent (or 'encode') the function by means of its graph. Representing the function $f_H$ as a set would produce

$$f_H = \{(f_H, 0)\}.$$

Recall that in $ZF$, the concept of an ordered pair can be defined by $(a,b) = \{\{a\}, \{a,b\}\}$, so that the ordered pair is a set, and $(a,b) = (c,d)$ if and only if

$a = c$ and $b = d$ (see [V, p. 15]). Then, the ordered pair $(f_H, 0)$ can be represented thus:

$$(f_H, 0) = \{\{f_H\}, \{f_H, 0\}\}$$

which would lead to:

$$f_H = \{(f_H, 0)\} = \{\{\{f_H\}, \{f_H, 0\}\}\}.$$

Let us denote: $A = \{\{f_H\}, \{f_H, 0\}\}$, $B = \{f_H\}$, and observe that $f_H = \{A\}$. Then we have

$$f_H \in B \in A \in f_H,$$

which would imply the following infinite $\in$-decreasing sequence

$$\cdots \in f_H \in B \in A \in f_H \in B \in A,$$

which stands in contradiction to the axiom of foundation.

In other words, by accepting $H$ as a legitimate proposition, we have arrived at a result that contradicts the axiom of foundation. Moreover, any proposition that is (strongly) self-referential would yield a similar result.

In conclusion, we have achieved the following result:

**Lemma 1:** Accepting a (strong) self-referential proposition as a legitimate sentence contradicts ZF's axiom of foundation (in the meta-model $N^*$).

From this we conclude: if we accept ZF's axiom of foundation and if we wish to determine that $N^*$ is consistent, we must disqualify self-referential propositions as meaningful propositions (that can be given a value of true or false) in $N^*$.

Let us add a few remarks regarding the constructions we performed in order to arrive at Lemma 1. We consider the following items.

1. Constructing the sets of type $P$

2. Constructing the set $F_P$.

3. Presenting every function from $F_P$ as a set.

First, the construction of type $P$ sets is natural and intuitive. There is a chance of circularity in the definition of these sets, which might lead these sets to be empty. However, assuming the existence of self-referential propositions like proposition $H$, which declares its own provability, ensures that there will be non-empty sets like these. More explicitly, accepting proposition $H$ as legitimate necessarily entails the existence of (non-empty) sets of type $P$.

Second, the construction of the set $F_P$, where the interpretation of a proposition is represented by a function, is also natural and analogous to what is common in computer programs, for example, where the semantics of the program is represented by a (full or partial) function $p: N \to N$ ($N$ here represents the set of natural numbers). We further recall that the assignment of truth value to sentences is represented by a function that assigns to each sentence a value from the set $\{0,1\}$. Recall also that the halting problem (in which the halting of a program is also an extensional property) can be presented as follows: 'does there exist a function $h(p,d)$, such that for every program $p$ and its input $d$, $h(p,d) = 1$ if $p$ halts at $d$ and $h(p,d) = 0$ if $p$ does not halt at $d$.'

Third, presenting any function from $F_P$ as a set is a common representation: every function is represented by a set of ordered pairs (known as the graph of the function), in which each ordered pair is itself a set. Furthermore, this representation is commonly used as the basic definition of a function (see [V, p. 16]).

Finally, when we combine the three stages of construction we find that, based on the assumption that the self-referential proposition $H$ is legitimate, we can define a set the existence of which contradicts the axiom of foundation, and the result is Lemma 1.

One may also claim that the above constructions result in an encoding of some special propositions by sets. We can argue that this coding is analogous to Gödel's coding, in which every formula in a first order language was coded as a natural number, or to a situation in which a computer program is coded by a natural number. This is an interesting analogy, but accepting it is not necessary to prove Lemma 1. It is enough to claim that from the existence of the non-empty set $A_P$, we can deduce the existence of the non-empty set $F_P$, and then the conclusion is drawn naturally from the common representation of functions as a set of ordered pairs.

**<u>Additional notes:</u>**

**Note 1:** The result presented in Lemma 1 can be extended in the following manner:

**<u>Lemma 2:</u>** In the meta-model $N^*$, every set of type $P$ is empty.

The proof for Lemma 2 is achieved similarly to the proof of Lemma 1, so we will content ourselves here with an abbreviated description.

If $A_P$ is a set of type $P$, then:

1. From Lemma 1, we deduce that there are no self-referential propositions in $A_P$.

2. Since every proposition in $A_P$ refers to the provability/non-provability of a different proposition, we must conclude that, assuming that $A_P$ is not empty, $A_P$ must include a proposition such as: $H_1$: The proposition $H_2$ is provable/unprovable. Since we must have $H_2 \in A_P$, we conclude that $H_2$ is a proposition such as: $H_2$: The proposition $H_3$ is provable/unprovable, where $H_3 \in A_P$.

3. If we continue in this manner, we will arrive at an infinite sequence of propositions of the form: $H_i$: the proposition $H_{i+1}$ is provable/unprovable (certain propositions may appear in the sequence more than once). Then, using explanations similar to those we used to prove Lemma 1, we will have an infinite sequence of functions, each contained in the domain of its predecessor. This will result in an infinite $\in$-decreasing sequence of sets, contradicting the axiom of foundation, and the result follows.

Lemma 2 leads us to the conclusion that, in addition to self-referential propositions, the existence of a (finite or infinite) sequence of propositions $H_1, H_2, \ldots, H_n, \ldots$ in which every proposition only refers to a given extensional property of a proposition in the sequence, also contradicts the axiom of foundation.

**Note 2:** It should be emphasized that the above results were obtained in the meta-model $N^*$. These results may not be valid in other models for which the meta-interpretation we used is not allowed.

**Note 3:** The results of Lemma 1 and Lemma 2 are significant only if we accept the axiom of foundation. There are areas in which it is reasonable not to rely on this axiom absolutely (for example in the semantics of natural language). On the other hand, in areas that require rigorousness, we are likely to rely on the axiom of foundation (we write 'likely' because there may be situations in mathematics where this axiom is not relied upon as in PA (Peano Arithmetic), for example). And yet, the axiom of foundation is generally accepted by those who engage in number theory as well, even though in many cases it is not needed.

**Note 4:** It should be emphasized that the results of this chapter imply that some propositions which were not considered paradoxical (like the propositions 'L: L is not provable' or 'K: K is true') actually contradict the axiom of foundation in $N^*$. If we wish to understand the semantic roots of this contradiction, we should attempt to follow the meaning of such sentences, an attempt which will lead us to a process of infinite regress. The discussion in the next chapter may help to clarify this issue.

**Note 5:** In continuation of the previous note, we can add that the fact that there is an analogy between the axiom of foundation and the concept of 'ungroundedness,' as it is commonly used by philosophers, has been noted before (see e.g., [He]). What we have done in this section can be seen as a rigorous discussion of that analogy. In other words, our discussion can help clarify the difference between rigorous systems that assume the axiom of foundation and other systems.

### 5) A discussion on the sentence $G_T$

In the previous section, we concluded that the existence of a (strong) self-referential proposition in the meta-model $N^*$ turns $N^*$ into an inconsistent model. Before moving on to a discussion of this conclusion's implications, we first devote this section to the connection between self-referential propositions, as presented in the previous section, and the sentence $G_T$ that was presented by Gödel in the proof of his first incompleteness theorem.

We recall Gödel's saying on the interpretation of the sentence $G_T$: 'We therefore have before us a proposition that says about itself that it is not provable.' On the other hand, there are various claims that $G_T$ is not self-referential, but these claims are generally considered not correct. A detailed discussion of this issue, which is accessible to a wider audience, can be found in [F, p. 44-46 and p. 84 -87]. On page 45 the author notes, regarding sentences built by the fixpoint construction (like $G_T$): '…the sentences constructed in the proof that every arithmetical property $P$ has a provable fixpoint are self-referential in a stronger sense.' We will not follow this detailed explanation here, but refer the reader to [F] for further details.

Gödel also talked of the relation of his sentence to the liar paradox that results from the sentence which says about itself that it is untrue. Others go even further and say: 'the Gödel sentence in effect results from this paradoxical sentence on substituting 'provable' for 'true'.' (see for example, [BBJ, p. 227]). So one may say that the sentence $G_T$ is equivalent to the sentence:

$L: L$ is unprovable.

However, this may lead to confusion, because though $G_T$ is self-referential, there is a substantial difference between the sentence $G_T$ and the proposition $L$.

Let us discuss the conceptual difference between these two, according to the discussion in [F, p. 84-87]. On page 85 the author says of 'non-formal' sentences like $L$: 'Such self-reference has been thought to be suspect for various reasons. In particular, the charge of leading to an infinite regress when one attempts to understand the statement has been leveled against it.' This means that if we look for the 'content' or the 'meaning' of the proposition, we will find ourselves in a 'process of infinite regress.' On the other hand, Gödel's sentence 'refers to a syntactic object, a sequence of symbols, and there is no infinite regress involved in establishing the reference on the phrase. '

According to the above observation, one may ask whether there is any connection between the discussion and conclusions of the previous sections and Gödel's sentence $G_T$. We claim that the answer to this question depends on the interpretation (or meta-interpretation) one ascribes to $G_T$. In order to discuss this issue further, we shall consider the sentence $G_T$ in more detail, and follow the discussion in Chapter 16 of [S]. We begin by recalling some definitions given there.

1. The *diagonalization* of a sentence $\varphi = \varphi(y)$ is $\exists y(y = \ulcorner \varphi \urcorner \wedge \varphi)$, where $\ulcorner \varphi \urcorner$ here stands in for the numeral for $\varphi$'s Gödel number – the diagonalization of $\varphi(y)$ is thus equivalent to $\varphi(\ulcorner \varphi \urcorner)$.

2. The relation $Gdl(m,n)$ holds just when $m$ is the super Gödel number for a PA proof of the diagonalization of the wff (well-formed formula) with Gödel number $n$.

3. $\text{Gdl}(x, y)$ stands for a $\Sigma_1$ wff which canonically captures $Gdl$.

Now we follow Gödel in constructing the corresponding wff

$U(y) =_{\text{def}} \forall x \neg \text{Gdl}(x, y)$.

For convenience, the above wff will be abbreviated by 'U' when we do not need to express that it contains 'y' free.

Finally, we diagonalize U itself, to give

$G = \exists y(y = \ulcorner U \urcorner \wedge U(y))$, which is equivalent to $U(\ulcorner U \urcorner)$.

This is the 'Gödel sentence' denoted by $G_T$ in this paper (we will follow the notations of [S] in the rest of this section). Unpacking that a bit we find that G is equivalent to

$\forall x \neg \text{Gdl}(x, \ulcorner U \urcorner)$.

We denote the latter equivalent wff also by G.

We give now a 'meta-interpretation' for G. In that meta-interpretation, we have that G 'says' that there is no PA proof for the diagonalization of U. But the diagonalization of U is G itself. So we may conclude that G 'says about itself that it is unprovable' (in the following, when we write 'proof,' we mean to say a 'PA proof'). If at this point we refer only to the 'syntactic meaning' of the formula, interpreted as the 'diagonalization of U,' we observe no 'infinite regress' in the meta-interpretation of G.

On the other hand, when we consider meta-interpretation in $N^*$, we can consider the 'content' or 'meaning' of a sentence. In particular, when we observe in such a meta-interpretation that a sentence 'talks' about the provability (or unprovability) of another sentence, questions like 'what is the meaning of that sentence' are legitimate. Furthermore, if we found in this meta-interpretation that this sentence also 'talks' about the provability of a sentence, we can continue questioning in the same manner.

Hence, using the above meta-interpretation for G, we may write the 'content' or 'meaning' of G as follows:

G: there is no proof to 'the diagonalization of U. '

However, since the diagonalization of U is G itself, which says that 'there is no proof to the diagonalization of U', then, if instead of considering the 'diagonalization of U' only as a syntactic object, we continue to 'meta-interpret' the 'diagonalization of U,' we may write the 'meaning' of G as follows:

G: there is no proof to 'there is no proof to the diagonalization of U, '

and applying this argument again:

G: there is no proof to 'there is no proof to 'there is no proof to the diagonalization of U. '

When we continue to apply the above argument, we conclude that there exists an infinite regress in this meta-interpretation of G.

Note further that the claim 'there is no proof to the 'diagonalization of U' is the only claim made by G (according to the above meta-interpretation), and since - again by the above meta-interpretation - the 'diagonalization of U' is G itself, we have that the sentence G is strongly self-referential, according to the denotation in Section 2.

We shall now summarize our observations. We may ascribe two meta-interpretations to the sentence G: a 'formal interpretation' and a 'referential interpretation.' Firstly, in both interpretations we ask: 'what does G say? ' and the answer is 'G says that there is no proof to the diagonalization of U.' Then comes a second question, which is different for each interpretation. In the 'formal interpretation' we ask: 'what is 'the diagonalization of U?' and the answer is: 'the diagonalization of U is G itself.' Now, since G is just a syntactic object, we have that G 'says' that there is no formal proof for the sentence G, and the result of the first incompleteness theorem follows. In particular, no infinite regress is observed under this interpretation. In the second 'referential interpretation,' we ask again the first question about the diagonalization of U: 'what does the diagonalization of U say? ' and the answer is 'the diagonalization of U says that there is no proof to the diagonalization of U,' and by repeating this question further, an infinite regress results.

It should be emphasized that the existence of the second meta-interpretation has nothing to do with the first incompleteness theorem, which remains a valid theorem. On the other hand, combining the second interpretation for G with the results of the previous section shows that the very existence of the Gödel's sentence causes a contradiction to the axiom of foundation in the model $N^*$.

### 6) A note regarding the significance of the results

In the wake of Gödel's incompleteness theorems, the opinion took root that sentences that can be interpreted as 'self-referential' (like $G_T$) are possible in mathematics. In the previous sections, we showed that this need not be always the case. On the one hand, in PA (Peano's axiom system for natural numbers), which does not include the axiom of foundation, such sentences are still possible. The same may be said of other formal theories which include PA. On the other hand, based on the second meta-interpretation of such sentences, which was introduced in the previous chapter, we deduce that accepting these sentences as legitimate in the meta-model $N^*$ may lead to a paradox. It follows that we should consider the possibility of rejecting the existence of (strong) self-referential sentences in the meta-model $N^*$ as natural and logical.

It is possible that further philosophical discussions on what we have said above are in place, but it is not our purpose to proceed with such discussions here. We only note that since the second interpretation introduced in the previous chapter uses arguments which are also used in the proof of the first incompleteness theorem, our opinion is that this second interpretation should not be ignored. From now on we limit ourselves to the following summary of the results up to this point:

1. Accepting a (strong) self-referential proposition in the meta-model $N^*$ as legitimate turns $N^*$ into an inconsistent model.

2. Under some meta-interpretation of Gödel's sentence $G_T$, a contradiction with the axiom of foundation occurs in the meta-model $N^*$.

These two conclusions will be our assumptions in the following two sections, where we present methods for addressing the problem that they raise. These methods can also be adapted to exclusively address the problem raised by result #1 (if we were to reject the interpretation upon which the second result relies). They may also be relevant if we reject both results, but are interested in methods for addressing the problem of the existence of 'strange' self-referential sentences (like $G_T$) in mathematics.

### 7) Approaches for coping with the problem

There are various possibilities for dealing with the problem of the existence of self-referential proposition in $N^*$. In what follows we present several suggestions that look worthy of discussion. They are:

1. Accepting the fact that $N^*$ contains 'paradoxical' propositions.

2. Renouncing the axiom of foundation in $N^*$.

3. Making a change in the first order language or in the theory $T$ that leads to the rejection of self-referential sentences, such as $G_T$.

4. Imposing limits on the interpretation (meta-interpretation) of sentences of $T$.

5. Accepting that certain sentences in the language will not have truth values (false or true). In other words, renouncing, to some extent, the law of excluded middle, and accepting that self-referential sentences will have a third value, designated, for example, NM (No Meaning).

In the following, we confine ourselves to the fifth possibility. The fact that this possibility waives the law of excluded middle seems like an extreme step. However, in light of the fact that there are well-formed sentences that lead to paradox (at least under a certain interpretation), the idea of adopting three truth values (True, False, NM) may be acceptable. This bitter pill can be sweetened with the clarification that, for all 'regular' sentences, the truth value would remain unchanged (i.e., True or False, as before). The new value (NM) would be reserved only for very special sentences (that are interpreted as self-referential).

Three-valued logic has been used more than once in addressing problems analogous to the problem that concerns us here (see, for instance, [K] for the use of Kleene's strong three-valued logic, or [B] for Bochvar's internal three-valued logic, also known as Kleene's weak three valued logic). This proposal appears to return us to the situation in natural language wherein propositions can have vague meaning, or lack meaning at all – something that we wish to avoid. In practice, however, if we define a logic that embraces the idea behind Bochvar's internal three-valued logic, we will find that the impact of using three-valued logic is not so revolutionary. Under such a logic, the discussion can effectively be reduced to sentences that are not interpreted as self-referential, for which the law of excluded middle remains very much alive.

Let us move on to describe the proposal itself. Accepting this proposal requires that we replace the law of excluded middle with a different law, which we call 'the extended law of excluded middle':

**The extended law of excluded middle:** For every sentence $A$ in a language, one (and only one) of the following is true:

1. The truth value of $A$ is NM.

2. $A \vee \neg A$.

In this logic, we demand that the truth value of a sentence that contains a component with the truth value NM will always be NM. Moreover, let us declare that the domain of the quantifiers will include only the values of the variables for which the truth value of the relevant sentence is other than NM.

The resulting state requires that we redefine what constitutes a consistent theory, and when a theory is complete. The definitions are as follows:

**Definition:** A theory is *consistent* if for every sentence $A$ there is no proof for both $A$ and $\neg A$, and if every sentence with a truth value NM has no proof.

**Definition:** A consistent theory is *complete* if for every sentence $A$ with a truth value that is not NM, there is proof for $A$ or $\neg A$.

As long as we wish that the meta model $N^*$ will be consistent, the following notes should be considered.

1. In a consistent theory, the truth value of every (strong) self-referential sentence must be NM.

2. Since the ('normal') rules of inference are valid only when they or their components have a truth value of True or False, we determine that a sentence with a truth value of NM has no proof, and it cannot appear as a component of any rule of inference in any consistent theory.

3. Therefore, in a consistent theory, a sentence with a truth value of NM cannot appear in any proof.

4. Moreover, in a consistent theory, a sentence with a truth value of NM cannot be used as an axiom (because axioms are automatically associated with the truth value True).

In a language that is informal, we can say, in light of the definitions and the comments we have presented, that sentences with a truth value of NM are 'excluded' in any consistent theory: they cannot be proved and nothing can be inferred from them. Therefore, at least in practical terms, we can make the rather vague argument that, under the new definitions, all of the mathematical work within the meta-model $N^*$ 'remains as it was.' The only difference between the new situation and the previous one is a newfound 'disregarding of self-referential sentences.'

### 8) The implications of the problem

The problem we have raised carries additional implications, which will be briefly discussed in this section.

**Firstly, it raises a question regarding the first incompleteness theorem itself:** What is the value of the theorem if we have shown, at least apparently, that there is a flaw in its proof? As we noted earlier, the answer to this question is simple: there is no flaw in the construction of the sentence $G_T$ and in the inference that it 'refers to itself and has no proof in $T$.' On the other hand, the theorem gives rise to a claim that is no less important than its original one: 'there exists a well-defined sentence in $T$ that is paradoxical, when meta-interpreted in $N^*$.'

We therefore find that under a certain interpretation the theorem faces us with a problem that is 'worse' than incompleteness, namely, inconsistency. As we have shown in the previous section, this problem can probably be circumvented by introducing changes to the theory or its interpretation. However, it will be necessary to check whether such circumvention will be able to prevent the phenomenon of self-reference in an entirely formal way (without the need address the sentences' meaning). If we adopt the suggestion of extending the law of excluded middle (as we did in the previous section), the conclusion of Gödel's theorem can be presented in the following manner (when $T$ is a theory which contains a 'sufficient part' of ZF - including the axiom of foundation):

If $T$ is consistent, then there exists a sentence that has 'no meaning' (NM) in $T$.

In conclusion, we find that the validity and the importance of Gödel's theorem remain unchanged, so long as we provide the theorem with a suitable interpretation. We might even claim that there is an interpretation which carries a 'more significant impact' – an impact that raises concerns regarding possible changes to the theory and/or its interpretation.

**Secondly, the problem returns us once again to the issue of incompleteness.** For example, if we adopt the definition of completeness from the previous section, we are once again faced with the question of whether a theory $T$ is complete (Gödel's theorem will not help us when the sentence $G_T$ has a truth value of NM, so its improvability does not lead to incompleteness under the new definition). However, there are proofs for incompleteness theorems that do not rely on self-reference. One such is, for instance, Chaitin's incompleteness theorem (see [Cha], or [F, chapter 8]), which relies on principles from Information Theory. Chaitin himself writes in [Cha] about the difference between his approach and Gödel's:

> Gödel's theorem may be demonstrated using arguments having an information-theoretic flavour. In such an approach it is possible to argue that if a theorem contains more information than a given set of axioms, then it is impossible for the theorem to be derived from the axioms. In contrast with the traditional proof based on the paradox of the liar, this new viewpoint suggests that the incompleteness phenomenon discovered by Gödel is natural and widespread rather than pathological and unusual.

### 9) Summary and additional notes

In this paper we proposed an interpretation for self-referential sentences in a 'meta-model' $N^*$ of ZF. Under this interpretation, we obtained the following two results:

1. Accepting (strong) self-referential sentences as legitimate turns the meta-model $N^*$ into an inconsistent model.

2. Under a certain interpretation of Gödel's sentence $G_T$, it is possible to argue that accepting this sentence as legitimate turns the meta-model $N^*$ into an inconsistent model.

These two conclusions were followed by a basic discussion (Section 7) which offered suggestions for addressing the problems that arise from the conclusions. We considered using the ideas of Bochvar's internal three-valued logic, which in general terms can be said to 'exclude self-referential sentences' and therefore 'solve the problem by "ignoring" it.' This approach is one of several different possibilities for managing the problem of self-reference (see for instance, [D, p. 10]), but it seems best suited for use in rigorous mathematical systems.

If we accept the interpretation introduced in this article, the above results raise additional questions and problems that must be addressed, but the discussion of which is beyond the scope of this paper. For example:

1. We must keep searching and inquiring into the possibilities for addressing the existence of paradoxical sentences in $N^*$, either in continuation of the options presented in Section 7 or in new and original directions. One specific question that must be addressed is whether it is possible to define a theory of arithmetic that does not allow for the construction of self-referential sentences.

2. It would be interesting to continue examining the influence of our discussion on the various existing proofs of incompleteness, and on the results and conclusions that would be produced by relying on those theorems.

3. Problems of self-reference appear elsewhere in mathematics. Some proofs that use fixpoint construction or other diagonal methods can be said to be involved with self-reference. It would therefore be interesting to study how the discussion in this paper might influence these proofs. Two important problems in this area are the halting problem in the Theory of Computability, and the use of Cantor's diagonal method to show that the set of real numbers is not countable. It is likely that the issues we have raised here will have an impact on these proofs. These topics are addressed in [L1] and [L2]. The discussion in this paper is relevant to self-reference problems in other (primarily philosophical) contexts as well, which have not been addressed here (see, for instance, [D], [He]).

**Arieh Lev** School of Computer Sciences, The Academic College of Tel-Aviv-Yafo, Tel-Aviv, Israel. arieh@mta.ac.il